\documentclass[11pt]{amsart}
 
\usepackage[arrow,matrix,curve,cmtip,ps]{xy}
\usepackage{upref}
\usepackage{fullpage}
\newtheorem{thm}{Theorem}[section]   
\newtheorem{lem}[thm]{Lemma}         
\newtheorem{prop}[thm]{Proposition}  

\theoremstyle{definition}
\newtheorem{defn}[thm]{Definition}   

\theoremstyle{remark}
\newtheorem{exs}[thm]{Examples}        

\numberwithin{equation}{section}     

\newcommand{\thmref}[1]{Theorem~\textup{\ref{#1}}}

\newcommand{\lemref}[1]{Lemma~\textup{\ref{#1}}}
\newcommand{\propref}[1]{Proposition~\textup{\ref{#1}}}

\renewcommand{\c}[1]{\mathcal #1}

\newcommand{\B}{\mathcal B}
\newcommand{\K}{\mathcal K}

\newcommand{\A}{\mathcal A}

\newcommand{\C}{\mathbb C}

\renewcommand{\phi}{\varphi}
\renewcommand{\theta}{\vartheta}

\newcommand{\midtext}[1]{\quad\text{#1}\quad}

\DeclareMathOperator{\Ad}{Ad}
\DeclareMathOperator{\ad}{Ad}
\DeclareMathOperator{\id}{id}

\DeclareMathOperator{\aut}{Aut}

\DeclareMathOperator{\UM}{{\mathit{UM}}}

\newcommand{\nonzero}[1]{\qquad\text{if $#1$
 \textup(and $0$ otherwise\textup)}}

\newcommand{\rip}[1]{_{\!#1}}

\newcommand{\<}{\langle}
\renewcommand{\>}{\rangle}

\newcommand{\what}{\widehat}
\renewcommand{\hat}{\widehat}

\renewcommand{\:}{\colon}
\newcommand{\deltahat}{{\widehat\delta}}
\newcommand{\epsilonhat}{{\widehat\epsilon}}
\newcommand{\deltahathat}{{\widehat{\widehat\delta}}}
\newcommand{\betahathat}{{\widehat{\widehat\beta}}}
\newcommand{\cotimes}{{\otimes_*}}

\begin{document}

\title{Maximal Coactions}

\author{Siegfried Echterhoff}
\address{Westf\"alische Wilhelms-Universit\"at\\
Mathmatisches Institut\\
Einsteinstr. 62\\
D-48149 M\"unster\\
Germany}
\email{echters@math.uni-muenster.de}

\author{S.~Kaliszewski}
\address{Department of Mathematics\\
Arizona State University\\
Tempe, Arizona 85287}
\email{kaliszewski@asu.edu}

\author{John Quigg}
\address{Department of Mathematics\\
Arizona State University\\
Tempe, Arizona 85287}
\email{quigg@math.la.asu.edu}


\date{September 12, 2001}

\begin{abstract}
A coaction $\delta$ of a locally compact group $G$ on a $C^*$-algebra
$A$ is \emph{maximal} if a certain natural map from $A\times_\delta
G\times_{\hat\delta}G$ onto $A\otimes\K(L^2(G))$ is an isomorphism.
All dual coactions on full crossed products by group actions
are maximal; a discrete coaction is maximal if and
only if $A$ is the full cross-sectional algebra of the corresponding
Fell bundle.  For every nondegenerate coaction of $G$ on $A$, there is
a maximal coaction of $G$ on an extension of $A$ such that the quotient
map induces an isomorphism of the crossed products.
\end{abstract}

\maketitle


\section{Introduction} 

One interesting thing about Katayama's duality theorem for coactions,
and more generally Mansfield's imprimitivity theorem, is that it 
appears at first to fundamentally involve \emph{reduced} 
crossed products by dual actions.  
This remains largely true even when the coactions are full 
(\emph{full} coactions are
defined using full group $C^*$-algebras; the original,
spatially-defined coactions are now called \emph{reduced}): 
in \cite{KQ-IC}, the second two authors 
proved that for any nondegenerate \emph{normal} 
full coaction
$(A,G,\delta)$ and any closed normal subgroup $N$ of $G$, the reduced
crossed product $A\times_{\delta}G\times_{\what\delta,r}N$ is
Rieffel-Morita equivalent to $A\times_{\delta|}G/N$.  
On the other hand, the first two authors, together with Iain~Raeburn,
showed in \cite{EKR-CP} that 
for any \emph{dual} coaction $\delta$, 
the \emph{full} crossed product
$A\times_{\delta}G\times_{\what\delta|}N$ is naturally 
Rieffel-Morita
equivalent to $A\times_{\delta|}G/N$.
An investigation of this phenomenon for
arbitrary coactions of discrete groups led the first and third authors
in \cite{EQ-FD} to conclude that in general,
Mansfield imprimitivity holds for neither the full nor the reduced
crossed product, but for an ``intermediate'' crossed product
$A\times_\delta G\times_{\hat\delta,\mu}N$.  
While that work relied heavily on the connection between discrete
coactions and Fell bundles, it also revealed a general framework which 
provided the immediate impetus for the work which appears here.
In particular, it motivated the definition and study of what we call
\emph{maximal} coactions.

Loosely speaking, the maximal coactions are those for which
full-crossed product duality holds; more precisely, $(A,G,\delta)$ is
maximal if a certain natural map from $A\times_\delta
G\times_{\what\delta}G$ onto $A\otimes\K(L^2(G))$ is an isomorphism.  
Thus, the maximal coactions are those for which the universal
properties of the full crossed product and the power of duality theory
can simultaneously be engaged.  
Our main theorem (\thmref{max thm}) says that for
every nondegenerate coaction $(A,G,\delta)$,
there is a maximal coaction on an extension of $A$ which has the same
crossed product;  we call this a \emph{maximalization} of $\delta$, and
show that it is unique up to an appropriate notion of isomorphism.  
Since there is always a coaction on a quotient of $A$ (the
\emph{normalization} of $A$) which has the same crossed product and for
which reduced-crossed-product duality holds, 
we can conclude that in general 
(as for the discrete coactions in \cite{EQ-FD}), Katayama duality
holds for
some intermediate crossed product lying between the full and reduced
crossed products by the dual action. 

The paper is organized as follows: in Section~2, we prove a lemma on
the naturality of normalizations
which, surprisingly, forms an essential part of the
proof of our main theorem.  We also show that the normal coactions are
\emph{precisely} those for which reduced-crossed-product duality holds. 
In Section~3, we define maximality and maximalization, and prove our
main theorem.  The basic ingredients of the proof are the observations
that dual coactions are always maximal (as suggested by the results in
\cite{EKR-CP}), and Rieffel-Morita equivalence preserves maximality.
Section~4 details the natural
relationship between maximalizations of
discrete coactions and the 
maximal cross-sectional algebras of the corresponding Fell bundles.

\subsubsection*{Acknowledgements}
Parts of this research were accomplished during separate visits of the
second and third authors to the first.
The second two authors are grateful to their hosts, 
and especially to Joachim Cuntz, for their generous
hospitality during their respective visits.
This research was partially supported by 
Deutsche Forschungsgemeinschaft (SFB 478).

 
\subsubsection*{Notation}

We adopt the conventions of \cite{EQ-IC, 
QuiFR,  
QuiDC} for coactions of groups on $C^*$-algebras, and of
\cite{ExeAF, FD-RA2} for Fell bundles. 
In particular, our coactions $(A,G,\delta)$ 
are \emph{full} in the sense that they map $A$
into $M(A\otimes C^*(G))$, the multiplier algebra of the minimal tensor
product of $A$ with the full group $C^*$-algebra.  Also, our coactions are
assumed to be \emph{nondegenerate}, in the usual sense that
the closure of
$\delta(A)(1\otimes C^*(G))$ is $A\otimes C^*(G)$.  
(When $B$ is an algebra and $X$ is a $B$-module, $BX$ denotes
the linear span of the set $\{ bx\mid b\in B, x\in X\}$.)
We use $(j_A,j_{C(G)})$ to denote the canonical covariant homomorphism
of $(A,C_0(G))$ into $M(A\times_\delta G)$, and $(i_B,i_G)$ for the
canonical covariant homomorphism of $(B,G)$ into $M(B\times_\alpha G)$
for an action $\alpha$ of $G$ on $B$. 
Additionally, we set $k_A = i_{A\times G}\circ j_A$, $k_{C(G)} =
i_{A\times G}\circ j_{C(G)}$, and $k_G = i_G$, 
which are maps of $A$, $C_0(G)$,
and $G$, respectively, into $M(A\times_{\delta}G\times_{\deltahat}G)$.  
We use superscripts, as in $k_G^A$, to distinguish between maps when
coactions on more than one algebra are commingled. 

If $\phi\colon A\to B$ is a homomorphism and $p\in M(A)$ is a
projection, we write $\phi|_p$ for the restriction of $\phi$ to the
corner $pAp$.  If $\delta$ is a coaction of $G$ on $A$ and $p$ is
\emph{$\delta$-invariant} in the sense that $\delta(p) = p\otimes 1$,
then we view $\delta|_p$ as a coaction on $pAp$; $\delta|_p$ will
be nondegenerate when $\delta$ is.

For a locally compact group $G$, we denote by $\lambda$ and $\rho$ the
left and right regular representations on $L^2(G)$, and by $\tau$ and
$\sigma$ the actions of $G$ on $C_0(G)$ by left and right translation,
respectively. $M:C_0(G)\to \B(L^2(G))$ denotes the representation
by multiplication operators,
and we use $u$ for the canonical map of $G$ into
$\UM(C^*(G))$.
By $\Sigma$ we mean the flip map: $A\otimes B\to B\otimes A$ for any
$A$ and $B$.

 
\section{Normal coactions} 
\label{normal} 
 
Let $\delta \: A \to M(A \otimes C^*(G))$ be a  
coaction of a locally compact group $G$ on $A$. 
Following \cite{NilDF} we define
a $^*$-homomorphism
$$\Phi=\pi\times U\colon A\times_{\delta}G\times_{\widehat{\delta}}G\to 
A\otimes
\K(L^2(G)),$$
where the covariant homomorphism $(\pi,U)$ of the dual system
$\big(A\times_{\delta}G, G,\widehat{\delta}\big)$ is defined by
\begin{equation} 
\label{eq-cov} 
\pi = (\id_A \otimes \lambda) \circ \delta \times (1 \otimes M),
\quad
U = 1 \otimes \rho.
\end{equation} 
As observed in \cite[Corollary~2.6]{NilDF}, $\Phi$ is always 
surjective, essentially because
$\overline{M(C_0(G))\rho(C^*(G))}=\K(L^2(G))$. In what follows we shall 
refer to
$\Phi$ as the {\em canonical surjection} of 
$A\times_{\delta}G\times_{\widehat{\delta}}G$ onto 
$A\otimes\K(L^2(G))$.

Recall from \cite[Definition~2.1]{QuiFR} that a coaction $(A,G,\delta)$ is 
\emph{normal} if the canonical map 
$j_A \: A \to M(A \times_\delta G)$  is injective.  
If $(A,G,\delta)$ is not normal, 
there always exists a normal coaction $\delta^n$ of $G$ on
$A^n=A/\ker j_A$ such that the quotient map of $A$ onto $A^n$ is 
$\delta - \delta^n$ equivariant and
induces an isomorphism of $A\times_\delta G$ onto
$A^n\times_{\delta^n}G$ (\cite[Proposition~2.6]{QuiFR}). 
We will call this
coaction $\delta^n$ the \emph{canonical normalization} 
of $\delta$.
In general, we say that $(B,G,\epsilon)$ is a
\emph{normalization} of $(A,G,\delta)$ if $\epsilon$ is normal and
there exists a $\delta - \epsilon$ equivariant surjection of $A$ onto
$B$ which induces an isomorphism between $A\times_\delta G$ and
$B\times_\epsilon G$. 

The following lemma implies that any two normalizations of a given
coaction are isomorphic; we state it in somewhat greater generality for
use in the proof of our main theorem.

\begin{lem}\label{normal-lem}
Let $(A,G,\delta)$ and $(C,G,\eta)$ be coactions, and let $\theta\colon
A\to C$ be a $\delta - \eta$ equivariant homomorphism.  If
$(B,G,\epsilon)$ is a normalization of $\delta$ with associated
equivariant surjection $\psi\colon A\to B$, and if $(D,G,\zeta)$ is a
normalization of $\eta$ with associated equivariant surjection
$\omega\colon C\to D$, there exists a unique homomorphism $\chi\colon
B\to D$ such that the diagram
\[
\xymatrix{
{A} 
\ar[rr]^{\theta}
\ar[d]_{\psi}
& 
& 
{C}
\ar[d]^{\omega}
\\
{B} 
\ar[rr]^{\chi}
& 
& 
{D}
}
\]
commutes.  The homomorphism $\chi$ will be $\epsilon - \zeta$ equivariant;
$\chi$ will be an isomorphism if $\theta$ is.
\end{lem}

\begin{proof}
We first consider the case $B=A^n$ and $D=C^n$; we need to show
that $\ker\psi\subseteq \ker(\omega\circ\theta)$, and this
is the same as showing $\ker j_A \subseteq\ker(j_C\circ\theta)$.  But
since $\theta$ is equivariant, the pair $(j_C\circ\theta,j_{C(G)}^C)$
is covariant for $(A,G,\delta)$, whence $j_C\circ\theta$ factors
through $j_A$ by the universal property of the crossed product.
This proves the desired inclusion.  Equivariance follows from a routine
diagram chase and surjectivity of $\psi$. 

We next consider the case $B=A^n$, $C=A$, and $\theta=\id$.  
Since $(D,G,\zeta)$ is
assumed normal, we have $D^n=D$, so the case proved above provides the
equivariant homomorphism 
$\chi\colon A^n\to D$, which is surjective because $\omega$ is.  
Since $\omega$ is equivariant, we have 
\[
(\omega\times G)\circ j_A = j_D\circ\omega,
\]
with $\omega\times G$ bijective because $\zeta$ is a normalization of
$\delta$, and $j_D$ injective because $\zeta$ is normal.   
Thus $\ker\omega = \ker j_A = \ker\psi$,
which implies that $\chi$ is injective.

The general case is obtained by cobbling these special cases
together. If $\theta$ is an isomorphism, applying the above to
$\theta^{-1}$ (together with the surjectivity of $\psi$ and $\omega$)
shows that $\chi$ is an isomorphism.
\end{proof}

Katayama's duality theorem 
(\cite[Theorem 8]{KatTD}; see also \cite[Corollary 2.6]{NilDF})
shows that crossed-product duality holds
for any normal coaction. The following proposition provides a converse
to this result:  a coaction $(A,G,\delta)$ is normal
if and only if Katayama duality holds in the sense that a certain
natural map is an isomorphism of $A\otimes\K(L^2(G))$ onto
$A\times_\delta G\times_{\what\delta,r}G$. 

\begin{prop} 
\label{prop-kat-reduced} 
Let $(A,G,\delta)$ be a coaction, 
let $\Lambda\colon  A \times_{\delta} G \times_{\what\delta} G \to 
 A \times_{\delta} G \times_{\what\delta,r} G $ be the regular
representation, 
let $\Phi \: A \times_{\delta} G \times_{\what\delta} G \to A 
\otimes  
\c K(L^2(G))$ be the canonical surjection, 
and let $\psi \: A \to A^n$  
denote the quotient map.  
\begin{enumerate}
\item
There exists an isomorphism
$\Upsilon$ of 
$A\times_\delta G\times_{\deltahat,r}G$ 
onto $A^n\otimes\K$ 
such that $(\psi\otimes\id_{\K})\circ\Phi = \Upsilon\circ\Lambda$.
\item
$\delta$ is normal if and only if the surjection 
$$\Psi= \Upsilon^{-1}\circ(\psi\otimes\id)\colon 
A\otimes\K \to A\times_\delta G\times_{\deltahat,r}G$$
is an isomorphism.
\end{enumerate}
\end{prop} 

The following diagram illustrates the proposition:
\[
\xymatrix{
{A\times_\delta G\times_{\deltahat}G} 
\ar[d]_{\Lambda}
\ar[rr]^{\Phi}
&
& {A\otimes\K(L^2(G))} 
\ar[d]^{\psi\otimes\id_\K}
\ar[dll]_{\Psi}
\\
{A\times_\delta G\times_{\deltahat,r}G}
\ar[rr]_{\Upsilon}^{\cong}
&
& {A^n\otimes\K(L^2(G)).}
}
\]
Although it is possible to fashion a proof from the existing theory ---
Katayama duality (as in \cite[Corollary~2.6]{NilDF} for nondegenerate 
normal full coactions), 
properties of normalizations, and naturality of the
crossed product --- we furnish a self-contained proof here.  
Thus, our argument provides an independent proof of Katayama duality.

\begin{proof}
To prove part~(i), it suffices to show that
\begin{equation}
\label{eq-ker}
\ker\Lambda = \ker((\psi\otimes\id_{\K})\circ\Phi)
\end{equation}
in $A\times_\delta G\times_{\deltahat}G$. 
Let $W=M\otimes\id_G(w_G)\in \UM(\K(L^2(G))\otimes C^*(G))$,
where $w_G$ denotes the canonical map $s\mapsto u_s$ of $G$ into 
$\UM(C^*(G))$,
regarded as an element of $\UM(C_0(G)\otimes C^*(G))$.
We claim that
\begin{equation}
\label{eq-ad}
\ad (1 \otimes W^*) \circ (\id_A \otimes \Sigma) \circ 
(\delta \otimes \id_{\c K}) \circ  \Phi
= (\pi \otimes 1) \times (U \otimes u)
\end{equation}
as maps from $A\times_\delta G\times_{\hat\delta}G$ to
$M(A\otimes\K\otimes C^*(G))$, 
with $\pi$ and $U$ as in \eqref{eq-cov}.  
For this, first observe that 
\begin{equation} 
\label{eq-sigma} 
\Sigma \circ (\id_G \otimes \lambda) \circ \delta_G 
= (\lambda \otimes \id_G) \circ \delta_G 
\midtext{and} 
\ad W^* \circ (\lambda \otimes \id_G) \circ \delta_G 
= \lambda \otimes 1. 
\end{equation} 
To see the latter represent $C^*(G)$ faithfully on a Hilbert space 
$\c  
H$ and compute for 
$\xi \in L^2(G, \c H) \cong L^2(G) \otimes \c  H$
and $s,t\in G$: 
\begin{align*} 
&\bigl( W^* \bigl( (\lambda \otimes \id_G) \circ \delta_G(s) \bigr) 
W \xi \bigr)(t) 
= \bigl( W^* (\lambda_s \otimes u_s) W \xi \bigr)(t) 
= u_{t^{-1}} \bigl( (\lambda_s \otimes u_s) 
W \xi \bigr)(t) 
\\&\quad= u_{t^{-1}} u_s (W \xi)(s^{-1}t) 
= u_{t^{-1}s} u_{s^{-1}t} \xi(s^{-1}t) 
= \bigl( (\lambda_s \otimes 1) \xi \bigr)(t). 
\end{align*} 

Now we can verify that \eqref{eq-ad} holds on generators:
for $a \in A$ we  compute 
\begin{align*} 
&\ad (1 \otimes W^*) \circ (\id_A \otimes \Sigma) 
\circ (\delta \otimes \id_{\c K}) \circ\Phi(k_A(a)) 
\\&\quad=
\ad (1 \otimes W^*) \circ (\id_A \otimes \Sigma) 
\circ (\delta \otimes \id_{\c K}) 
\bigl( (\id_A \otimes \lambda) \circ \delta(a) \bigr) 
\\&\quad=
\ad (1 \otimes W^*) \circ (\id_A \otimes \Sigma) 
\circ (\delta \otimes \lambda) (\delta(a)) 
\\&\quad=
\ad (1 \otimes W^*) \circ (\id_A \otimes \Sigma) 
\circ (\id_A \otimes \id_G \otimes \lambda) 
\circ (\delta \otimes \id_G) (\delta(a)) 
\\&\quad=
\ad (1 \otimes W^*) \circ (\id_A \otimes \Sigma) 
\circ (\id_A \otimes \id_G \otimes \lambda) 
\circ (\id_A \otimes \delta_G) (\delta(a)) 
\\&\quad\overset{\makebox[0mm]{\scriptsize\eqref{eq-sigma}}}{=} 
\ad (1 \otimes W^*) \circ (\id_A \otimes \lambda \otimes \id_G) 
\circ (\id_A \otimes \delta_G) (\delta(a)) 
\\&\quad\overset{\makebox[0mm]{\scriptsize\eqref{eq-sigma}}}{=} 
(\id_A \otimes \lambda) (\delta(a)) \otimes 1 
\\&\quad=
(\pi\otimes 1)\times(U\otimes u)(k_A(a)),
\end{align*} 
while for $f \in C_0(G)$ we have 
\begin{align*} 
&\ad (1 \otimes W^*) \circ (\id_A \otimes \Sigma) 
\circ (\delta \otimes \id_{\c K}) \circ\Phi(k_{C(G)}(f)) 
\\&\quad= \ad (1 \otimes W^*) \circ (\id_A \otimes \Sigma) 
\circ (\delta \otimes \id_{\c K}) (1 \otimes M_f) 
\\&\quad= \ad (1 \otimes W^*) (1 \otimes M_f \otimes 1) 
\\&\quad= 1 \otimes M_f \otimes 1 
\\&\quad= (\pi\otimes 1)\times(U\otimes u)(k_{C(G)}(f)). 
\end{align*} 
For $s \in G$ and $\xi \in L^2(G, \c H)$ we have 
\begin{align*} 
&\bigl( W^* (\rho_s \otimes 1) W \xi \bigr)(t) 
= u_{t^{-1}} \bigl( (\rho_s \otimes 1) W \xi \bigr)(t) 
= \Delta(s)^{1/2} u_{t^{-1}} (W \xi)(ts) 
\\&\quad= \Delta(s)^{1/2} u_{t^{-1}} u_{ts} \xi(ts) 
= \bigl( (\rho_s \otimes u_s) \xi \bigr)(t), 
\end{align*} 
which implies 
\begin{align*} 
&\ad (1 \otimes W^*) \circ (\id_A \otimes \Sigma) 
\circ (\delta \otimes \id_{\c K})\circ\Phi(k_G(s)) 
\\&\quad= \ad (1 \otimes W^*) \circ (\id_A \otimes \Sigma) 
\circ (\delta \otimes \id_{\c K}) (1 \otimes \rho_s) 
\\&\quad= \ad (1 \otimes W^*) (1 \otimes \rho_s \otimes 1) 
\\&\quad= 1 \otimes \rho_s \otimes u_s 
\\&\quad= (\pi\otimes 1)\times(U\otimes u)(k_G(s)).
\end{align*} 

Now applying $\id_{A\otimes\K}\otimes\lambda$ to both
sides of \eqref{eq-ad} yields
\[ 
\Theta\circ \bigl( (\id_A \otimes \lambda) \circ \delta 
\otimes \id_{\c K} \bigr)\circ\Phi
=(\pi\otimes 1)\times(U\otimes\lambda),
\] 
where $\Theta$ is the invertible map
$\ad(1\otimes(\id_{\K}\otimes\lambda)(W^*))
\circ(\id_A\otimes\Sigma)$.
Since $(\pi \otimes 1) \times (U \otimes \lambda)$ is equivalent to  
the regular representation of $A \times_\delta G 
\times_{\what\delta}  
G$ induced from $\pi$, 
for \eqref{eq-ker} it is therefore enough to show that
$\ker((\id_A\otimes\lambda)\circ \delta \otimes\id_\K)
\circ \Phi  
= \ker (\psi \otimes \id_{\c  K}) \circ \Phi$.
But $\pi = (\id_A\otimes\lambda)\circ\delta\times(1\otimes M)$
is a faithful representation of $A\times_\delta G$, 
so we can identify the
canonical map $j_A\colon A\to M(A\times_\delta G)$ with
$(\id_A\otimes\lambda)\circ\delta$.  Thus, by definition, $\ker\psi =
\ker j_A = \ker (\id_A\otimes\lambda)\circ\delta$, 
so (since $\K$ is nuclear) we are done with part~(i).

Since $\delta$ is normal if and only if $\psi$ is an isomorphism, 
part~(ii) follows immediately from part~(i). 
\end{proof} 


\label{maximal}
\section{Maximal coactions}
 
\begin{defn}
\label{def-max}
Let $\delta\:A\to M(A\otimes C^*(G))$ be a coaction.
We say that $\delta$ is {\em maximal} if the canonical map
$\Phi\:A\times_{\delta}G\times_{\widehat{\delta}}G\to 
A\otimes\K(L^2(G))$
is an isomorphism.
A maximal 
coaction $(B,G,\epsilon)$ is a \emph{maximalization} of $\delta$ if
there exists an $\epsilon - \delta$ equivariant surjection
$\theta\colon B\to A$ such that the induced map $\theta\times G\colon
B\times_\epsilon G\to A\times_\delta G$ is an isomorphism. 
\end{defn}

It will follow from \thmref{max thm} that every coaction $(A,G,\delta)$
sits ``between'' a maximal coaction and a normal one, in the sense that
there exists a maximal coaction $(A^m,G,\delta^m)$ and a normal
coaction $(A^n,G,\delta^n)$, together with equivariant surjections:
$\phi\colon A^m\to A$ and $\psi\colon A\to A^n$.  

\begin{exs}
If $G$ is any locally compact group, $(C^*(G),G,\delta_G)$ is maximal, 
because $\delta_G$ is the dual coaction on $C^*(G) = \C\times_{\id}G$, 
and dual coactions are always maximal (\propref{prop-max}). 
The normalization of $\delta_G$ is the coaction $\delta_G^n$ of $G$ on
$C^*_r(G)$ determined
by $\lambda(s)\mapsto \lambda(s)\otimes u(s)$ 
(see \cite[Example~2.12]{QuiFC}).  Thus, if $G$
is nonamenable, $\delta_G$ is maximal but not
normal, and $\delta_G^n$ is normal but not maximal.  See the discussion
following \cite[Proposition~3.12]{QuiFC} for
a coaction which is neither maximal nor normal.
\end{exs}

\begin{thm}
\label{max thm}
Every nondegenerate
coaction $(A,G,\delta)$ has a maximalization.  
If $(B,G,\epsilon)$ and $(C,G,\eta)$ are two maximalizations
of $(A,G,\delta)$ with canonical equivariant surjections $\phi\colon
B\to A$ and $\theta\colon C\to A$, then there
exists a $\epsilon - \eta$ equivariant isomorphism $\chi$ of
$B$ onto $C$ such that $\theta\circ\chi= \phi$. 
\end{thm}

The idea of the proof is as follows: if $(A^m,G,\delta^m)$ were a
maximalization of $(A,G,\delta)$, then it would satisfy
$$A\times_\delta G\times_{\deltahat}G \cong 
A^m\times_{\delta^m}G\times_{\what{\delta^m}}G
\cong A^m\otimes \K(L^2(G)).$$
Thus $A^m$ would be retrievable from $A\times_\delta
G\times_{\deltahat}G$ by 
taking a rank-one projection $P$ in $M(\K)$ and then cutting down by
the image $p$ of $1\otimes P$ in $M(A\times_\delta
G\times_{\deltahat}G)$. 
A natural candidate for $\delta^m$ would then be
the restriction to $A^m$ of the double-dual coaction
$\deltahathat$ on $A\times_\delta G\times_{\deltahat}G$, but
a technicality arises here 
because in general $p$ won't be $\deltahathat$-invariant;
an adjustment by a one-cocycle makes up for
this defect.  

We begin by showing that dual coactions are always maximal.
(It was shown in \cite[Theorem~3.7]{QuiFC} that
$A\times_{\delta}G\times_{\what\delta}G\cong A\otimes\K(L^2(G))$ for
any dual coaction $\delta$, but the canonical map $\Phi$ was not
explicitly identified as the isomorphism.)

\begin{prop}
\label{prop-max}
Let $\beta \: G \to \aut B$ be an action.
Then the dual coaction 
$$\what\beta = (i_B  
\otimes 1) \times (i_G \otimes u) \: B\times_\beta G \to M(B  
\times_\beta G \otimes C^*(G))$$ 
on the full crossed product $B \times_\beta G$ 
is maximal.
\end{prop} 
\begin{proof} 
Let 
$ \Phi \:  
(B \times_\beta G) \times_{\what\beta} G  
\times_{\what{\what\beta}} G 
\to (B \times_\beta G) \otimes \c K(L^2(G)) $
be the canonical map for $\widehat\beta$.
Since $(\id_{B \times_\beta G} \otimes \lambda) \circ \what\beta =  
(i_B \otimes 1) \times (i_G \otimes \lambda)$, we have
\[ 
\Phi = (i_B \otimes 1) \times (i_G \otimes \lambda) 
\times(1 \otimes M) \times (1 \otimes \rho). 
\] 
We have to check that this map is injective. 
 
For this recall first that the Imai-Takai duality theorem 
(see \cite[Theorem~5.1]{RaeCR})
provides an isomorphism $\Omega\: B \times_\beta G 
\times_{\what\beta}  
G \to B \otimes \c K$ which is given as the integrated form  
\[ 
\Omega = ((\id_B \otimes M)\circ\tilde\beta)\times
(1 \otimes \lambda)\times (1 \otimes M) 
\] 
where for $b\in B$, 
$\tilde\beta(b) \in C^b(G,B)\subseteq M(B \otimes C_0(G))$ 
denotes  
the function $t \mapsto \beta_{t^{-1}}(b)$.  
Recall also that $\Omega$  
transports the double dual action $\what{\what\beta}$ to the 
diagonal  
action $\beta \otimes \ad \rho$ of $G$ on $B \otimes \c K$, 
so  we get an isomorphism 
\[ 
\Omega \times G \: B \times_\beta G \times_{\hat\beta} G 
\times_{\betahathat} G 
\to (B \otimes \c K) \times_{\beta \otimes \ad \rho} G. 
\] 
Since $1 \otimes \rho$ implements an exterior equivalence between  
$\beta \otimes \id$ and $\beta \otimes \ad \rho$, we also have an  
isomorphism 
\[ 
\Theta = (i_B \otimes \id_{\c K}) \times (i_G \otimes \rho) \: 
(B \otimes \c K) \times_{\beta \otimes \ad \rho} G 
\to (B \times_\beta G) \otimes \c K. 
\] 
The composition $\Theta \circ (\Omega \times G)$ is thus an isomorphism  
between 
$B \times_\beta G \times_{\hat\beta} G \times _{\hat{\hat\beta}} G$ 
and 
$(B \times_\beta G) \otimes \c  K$, 
and all we have to do is to check that this isomorphism  
has the same kernel as $\Phi$. 
 
In order to do this, we fix a faithful representation $\sigma\times V$ of
$B \times_\beta G$ on a 
Hilbert space $\c H$,
and define a unitary operator $R$ on $\c H \otimes L^2(G) \cong  
L^2(G,\c H)$ by $(R \xi)(t) = V_t \xi(t)$ for $\xi \in  
L^2(G,\c H)$.  We claim that 
\[
\ad R \circ (\sigma \times V \otimes \id_{\c K}) 
\circ \Theta \circ (\Omega\times G) 
= (\sigma \times V \otimes \id_{\c K}) 
\circ \Phi;
\]
for this it is enough to check that both maps do the same on the  
generators of 
$B \times_\beta G \times_{\hat\beta} G \times _{\hat{\hat\beta}} G$.
For example, if $\ell_B \:  
B \to M(B \times_\beta G \times_{\hat\beta} G \times _{\hat{\hat\beta}} G)$ denotes the canonical  
imbedding, for all $b \in B$ we have 
\begin{eqnarray*} 
\ad R \circ (\sigma \times V \otimes \id_{\c K}) \circ \Theta 
\circ (\Omega\times G)(\ell_B(b)) 
&=& \ad R \circ (\sigma \otimes \id_{\c K}) \circ \Omega(k_B(b))\\
&=&\ad R \circ (\sigma \otimes M)(\tilde\beta(b)) \\
&\overset{(*)}{=} &\sigma(b) \otimes 1 \\
&=&(\sigma \times V \otimes \id_{\c K}) ((i_B(b) \otimes 1)\\
&=& (\sigma\times V\otimes\id_{\K})\circ\Phi(\ell_B(b)),
\end{eqnarray*} 
where the starred equality follows from the calculation
\[ 
\bigl( R (\sigma \otimes M)(\tilde\beta(b)) R^* \xi \bigr)(t) 
= V_t \sigma(\beta_{t^{-1}}(b)) V^*_t \xi(t) 
= \sigma(b) \xi(t) 
\] 
for $b\in B$, $\xi \in L^2(G,\c H)$, and $t\in G$. 
We omit the easier computations on the other generators. 
\end{proof} 

Let us say that two coactions $(A,G,\delta)$ and $(B,G,\epsilon)$ are
{\em Morita equivalent} if there exists a coaction $(C,G,\eta)$ and 
$\eta$-invariant full projections $p,q\in M(C)$ such that 
$A=pCp$, $B=qCq$, 
$\eta|_A = \delta$, and $\eta|_B = \epsilon$.  It is not hard to see,
using a linking-algebra argument, 
that this definition agrees with those already in the literature
(see \cite{BS-CH,BuiME,ER-MI,KQR-DR,NgCC2}).

\begin{prop}
\label{ME}
Suppose $(A,G,\delta)$ and $(B,G,\epsilon)$ are Morita equivalent
coactions.  Then $\delta$ is maximal if and only if $\epsilon$ is.
\end{prop}

\begin{proof}
It suffices to consider the case where $p\in M(B)$ is an 
$\epsilon$-invariant full
projection, $A=pBp$, and $\delta=\epsilon|_A$; 
then the isomorphisms 
$A\times_\delta G\times_\deltahat G\cong 
k_B(p)(B\times_\epsilon G\times_\epsilonhat G)k_B(p)$ 
and $A\otimes\K\cong
(p\otimes1)(B\otimes\K)(p\otimes1)$ transport
the canonical surjection $\Phi_A$ to the
restriction of $\Phi_B$ to $A\times_\delta G\times_\deltahat G$.  
Thus $\ker\Phi_A$
is associated to $\ker\Phi_B$ under the Rieffel correspondence set up
by the full projection $k_B(p)$,  
so $\ker\Phi_A = \{0\}$ if and only if
$\ker\Phi_B=\{0\}$,
which completes the proof.
\end{proof}

Adapting
the definition from 
\cite[Definition~2.7]{LPRS-RC}, 
where it appears for reduced coactions, we 
say that a  unitary $U$ in $M(A\otimes C^*(G))$ is a \emph{1-cocycle}
for a coaction $(A,G,\delta)$ if 
\begin{enumerate}
\item
$\id\otimes\delta_G(U) = (U\otimes 1)(\delta\otimes\id(U))$, and
\item
$U\delta(A)U^*(1\otimes C^*(G)) \subseteq
A\otimes C^*(G)$.
\end{enumerate}
Two coactions $\delta$ and $\epsilon$ of $G$ on $A$ are \emph{exterior
equivalent} if there exists a 1-cocycle $U$ for $\delta$ such that
$\epsilon = (\ad U)\circ\delta$; in this case, $\delta$ is
nondegenerate if and only if $\epsilon$ is
(more generally, this is true for Morita equivalent coactions
--- see \cite[Proposition~2.3]{KQ-IC}).
Exterior equivalent coactions have naturally isomorphic crossed
products: to see this,
realize $A\times_\delta G$ and $A\times_\epsilon G$ as
subalgebras of $M(A\otimes\K(L^2(G))$, and then argue exactly as in the
proof of \cite[Theorem~2.9]{LPRS-RC} (where the same result is proved
for reduced coactions) that $\id_A\otimes\lambda(U)$
conjugates one to the other.

\begin{lem}
\label{V-lem}
$V = (k_{C(G)}\otimes\id)(w_G)$ is a 1-cocycle for $\deltahathat$; hence
$\tilde\delta = \Ad(V)\circ\deltahathat$ is a 
\textup(nondegenerate\textup) coaction
of $G$ on $A\times_\delta G\times_{\deltahat}G$. 
\end{lem}

\begin{proof}
To verify condition~(i), we have:
\begin{eqnarray*}
\id\otimes\delta_G(V) 
& = & \id\otimes\delta_G(k_{C(G)}\otimes\id(w_G))\\
& = & k_{C(G)}\otimes\id\otimes\id(\id\otimes\delta_G(w_G))\\
& = & k_{C(G)}\otimes\id\otimes\id((w_G)_{12}(w_G)_{13}),
\end{eqnarray*}
using the identity $\id\otimes\delta_G(w_G) = (w_G)_{12}(w_G)_{13}$.  
Now from the definitions,
$$\deltahathat\otimes\id(V) = ((k_A\otimes 1)\times (k_{C(G)}\otimes
1)\times (k_G\otimes u))\otimes\id(k_{C(G)}\otimes\id(w_G))
= k_{C(G)}\otimes\id\otimes\id((w_G)_{13}),$$
and clearly $V\otimes 1 = k_{C(G)}\otimes\id\otimes\id((w_G)_{12})$.

For condition~(ii), 
note that for $f\in C_0(G)$ and $z\in C^*(G)$, the product
$(f\otimes 1)w_G^*(1\otimes z)\in M(C_0(G)\otimes C^*(G))$ is actually
in $C_0(G)\otimes C^*(G)$, since it corresponds
to the function $s\mapsto f(s)u_s^*z$ in $C_0(G,C^*(G))$.  Thus
$(C_0(G)\otimes 1)w_G^*(1\otimes C^*(G))\subseteq C_0(G)\otimes
C^*(G)$.
Now temporarily let $B=A\times_\delta G\times_{\deltahat}G$ and compute:
\begin{eqnarray*}
\deltahathat(B)V^*(1\otimes C^*(G))
& = & \deltahathat(B)
\bigl(k_{C(G)}(C_0(G))\otimes 1\bigr)
k_{C(G)}\otimes\id(w_G^*)(1\otimes C^*(G))\\
& = & \deltahathat(B)k_{C(G)}\otimes\id
\bigl( (C_0(G)\otimes 1)w_G^*(1\otimes C^*(G))\bigr)\\
& \subseteq & \deltahathat(B)k_{C(G)}\otimes\id
(C_0(G)\otimes C^*(G))\\
& = & \deltahathat(B)(1\otimes C^*(G))\\
&\subseteq & B\otimes C^*(G).
\end{eqnarray*}

Nondegeneracy of $\tilde\delta$ follows from that of $\deltahathat$.
\end{proof}

\begin{lem}
\label{k-lem}
\begin{enumerate}
\item
The pair $(k_{C(G)},k_G)$ is a covariant homomorphism of
$(C_0(G),G,\sigma)$ into $M(A\times_\delta G\times_{\deltahat}G)$.
\item
$\tilde\delta(x) = x\otimes 1$ for all 
$x\in (k_{C(G)}\times k_G)(C_0(G)\times_\sigma G)$.
\end{enumerate}
\end{lem}

\begin{proof}
Assertion (i) is a straightforward consequence of the definition of
$\deltahat$.  
It suffices to check~(ii) on generators.  For $f\in C_0(G)$,
\begin{eqnarray*}
\tilde\delta(k_{C(G)}(f)) 
& = & V\deltahathat(k_{C(G)}(f))V^*\\
& = & k_{C(G)}\otimes\id(w_G)(k_{C(G)}(f)\otimes
1)k_{C(G)}\otimes\id(w_G^*)\\
& = & k_{C(G)}\otimes\id(w_G(f\otimes 1)w_G^*)\\
& = & k_{C(G)}(f)\otimes1,
\end{eqnarray*}
since $C_0(G)\otimes 1$ commutes with $M(C_0(G)\otimes C^*(G))$.
For $s\in G$ we have
\begin{eqnarray*}
\tilde\delta(k_G(s)) 
& = & k_{C(G)}\otimes\id(w_G)(k_G(s)\otimes
u(s))(k_{C(G)}\otimes\id(w_G^*))\\
& = & k_{C(G)}\otimes\id(w_G) k_G\otimes\id(\delta_G(s)) 
(k_{C(G)}\otimes\id(w_G^*))\\
&\overset{\dag}{=}& k_G(s)\otimes 1.
\end{eqnarray*}
The identity at $\dag$, which almost says that the pair $(k_G,k_{C(G)})$
is covariant for $(C^*(G),G,\delta_G)$, follows from 
part~(i) the same way the covariance of any homomorphism $(U,\mu)$ of
$(C^*(G),G,\delta_G)$
follows from that of $(\mu,U)$ for $(C_0(G),G,\tau)$  (see
\cite[Example~2.9(1)]{RaeCR}).  
\end{proof}

For brevity, let $\delta\cotimes \id$ denote the coaction
$(\id_A\otimes\Sigma)\circ(\delta\otimes\id_{\K})$  of $G$ on
$A\otimes\K$.  

\begin{lem}
\label{invt}
$\Phi$ is a $\tilde\delta - \delta\cotimes\id$ equivariant surjection
of $A\times_\delta G\times_{\deltahat}G$ onto $A\otimes\K$. 
\end{lem}

\begin{proof}
We only need to check equivariance; that is, 
we need to show that
$(\id_A\otimes\Sigma)\circ(\delta\otimes\id_K)\circ\Phi
=(\Phi\otimes\id)\circ\tilde\delta$.
But by \eqref{eq-ad}, we have
$(\id_A\otimes\Sigma)\circ(\delta\otimes\id_{\K})\circ\Phi
= \ad (1\otimes W)\circ ((\pi\otimes 1)\times (U\otimes u))$,
and straightforward calculations verify that
\begin{eqnarray*}
\ad (1\otimes W)\circ((\pi\otimes1)\times(U\otimes u))
& = & \ad (1\otimes W)\circ(\Phi\otimes\id)\circ\deltahathat\\
& = & (\Phi\otimes\id)\circ(\ad V)\circ\deltahathat\\
& = & (\Phi\otimes\id)\circ\tilde\delta.
\end{eqnarray*}
\end{proof}

\begin{proof}[Proof of \thmref{max thm}]
Fix, for the entire proof, a rank-one projection $P\in \K(L^2(G))$, 
let $q=(M\times\rho)^{-1}(P) \in C_0(G)\times_\sigma G$, and
let $p=(k_{C(G)}\times k_G)(q)$, which is a
$\tilde\delta$-invariant  projection in $M(A\times_\delta
G\times_{\deltahat}G)$ by \lemref{k-lem}.  
We may therefore define a nondegenerate coaction $(A^m,G,\delta^m)$ 
by setting
\[
A^m = p(A\times_\delta G\times_{\deltahat}G)p
\qquad{\rm and}\qquad \delta^m = \tilde\delta|_{p}.
\]
Notice that, 
by the definition of $\Phi$,  we have
\begin{equation}
\label{Phi-eq}
\Phi(p) 
=\Phi\circ(k_{C(G)}\times k_G)(q)
=1\otimes(M\times\rho)(q)
= 1\otimes P.
\end{equation}

Now $p$ is a full projection in 
$M(A\times_\delta G\times_{\deltahat}G)$, since if we put $C=
k_{C(G)}\times k_G(C_0(G)\times_\tau G)$, then clearly
$\overline{CpC}=C$, whence
\begin{align*}
\overline{(A\times_\delta G\times_{\hat\delta}G)p(A\times_\delta
G\times_{\hat\delta}G)}
&= 
\overline{(A\times_\delta G\times_{\hat\delta}G)CpC(A\times_\delta
G\times_{\hat\delta}G)}\\
&=
\overline{(A\times_\delta G\times_{\hat\delta}G)C(A\times_\delta
G\times_{\hat\delta}G)}
=A\times_\delta G\times_{\hat\delta}G.
\end{align*}
Thus
$\delta^m$ is Morita equivalent to
$\tilde\delta$,
which is in turn Morita equivalent
(in fact exterior equivalent) to $\deltahathat$, 
which is maximal by \propref{prop-max}.
It follows from \propref{ME} that $\delta^m$ is also maximal. 

Now, identifying $(1\otimes P)(A\otimes\K)(1\otimes
P)$ with $A$, it follows from \lemref{invt} that the restriction
$\Phi|_p$ of $\Phi$ to $A^m$ is a $\delta^m - \delta$ equivariant
surjection of $A^m$ onto $A$.  
To prove that $(A^m,\delta^m)$ is a maximalization of $(A,\delta)$,
it remains to show that
the integrated form $\Phi|_p\times G$ is an isomorphism of
$A^m\times_{\delta^m}G$ onto $A\times_\delta G$.  

For this, we first point out that if $\epsilon$ and 
$\ad U\circ\epsilon$ are exterior equivalent coactions of $G$ on $B$, and
$(C,G,\eta)$ is a normalization of $\epsilon$ with associated
surjection $\chi\colon B\to C$, then $(C,G,\ad V\circ\eta)$ is a
normalization of $\ad U\circ\epsilon$ with the same surjection, where
$V = \chi\otimes\id(U)$.  (This is a straightforward calculation.)
In particular, since the dual coaction $\deltahathat_r$ of $G$ on the
reduced crossed product $A\times_\delta G\times_{\deltahat,r}G$ is a
normalization of $\deltahathat$ 
with associated surjection $\Lambda$
(\cite[Propositions~2.3 and~2.6]{QuiFR}),
$(A\times_\delta G\times_{\deltahat,r}G,G,\eta)$
is a normalization of $(A\times_\delta
G\times_{\deltahat}G,G,\tilde\delta)$ with the same surjection, for the
appropriate choice of $\eta$.  

Moreover, it is easily seen that
$(A^n\otimes\K,G,\delta^n\cotimes \id)$ is a normalization of
$(A\otimes\K,G,\delta\cotimes \id)$, with surjection $\psi\otimes\id_{\K}$.  
It follows that the isomorphism $\Upsilon\colon A\times_\delta
G\times_{\deltahat,r}G\to A^n\otimes\K$ 
of \propref{prop-kat-reduced} is the unique homomorphism 
between normalizations provided
by \lemref{normal-lem},
and hence is 
$\eta - \delta^n\cotimes \id$ equivariant.

It now follows that
$\Phi\times G$ is an isomorphism, since by naturality
it is the composition of isomorphisms
\begin{multline*}
(A\times_\delta G\times_{\deltahat}G)\times_{\tilde\delta} G
\stackrel{\Lambda\times G}{\cong}
(A\times_\delta G\times_{\deltahat,r} G)\times_\eta G\\
\stackrel{\Upsilon\times G}{\cong} 
(A^n\otimes\K)\times_{\delta^n\cotimes \id}G
\stackrel{(\psi\otimes\id)\times G}{\cong} 
(A\otimes\K)\times_{\delta\cotimes \id}G.
\end{multline*}
Since  $\Phi\times G$ takes the image $\ell(p)$ of $p$ in $M(A\times_\delta
G\times_{\deltahat}G\times_{\tilde\delta}G)$ to 
$j_{A\otimes\K}(1\otimes P)\in
M((A\otimes\K)\times_{\delta\cotimes \id}G)$, naturality (again) implies
that $\Phi|_p\times G$ is the composition of isomorphisms
\begin{eqnarray*}
A^m\times_{\delta^m}G
&=&
p(A\times G\times G)p\times_{\tilde\delta|_{p}}G\\
&\cong& \ell(p)(A\times G\times G\times_{\tilde\delta} G)\ell(p)\\
&\stackrel{(\Phi\times G)|_{\ell(p)}}{\cong}&
j_{A\otimes\K}(1\otimes P)((A\otimes\K)
\times_{\delta\cotimes \id}G)j_{A\otimes\K}(1\otimes P)\\
&\cong& (1\otimes P)(A\otimes\K)
(1\otimes P)\times_{(\delta\cotimes \id)|_{1\otimes P}} G\\
&\cong& A\times_\delta G.
\end{eqnarray*}
This completes the proof that $\delta^m$ is a maximalization of
$\delta$.

For the last statement of the theorem, consider the diagram
\[
\xymatrix{
{B\times_\epsilon G\times_{\hat\epsilon}G,\tilde\epsilon} 
\ar[rr]^{\phi\times G\times G}
\ar[dd]_{\Phi_B}
&
&
{A\times_\delta G\times_{\hat\delta}G,\tilde\delta} 
\ar[d]^{\Phi_A}
&
&
{C\times_\eta G\times_{\hat\eta }G,\tilde\eta} 
\ar[ll]_{\theta\times G\times G}
\ar[dd]^{\Phi_C}
\\
& 
& 
{A\otimes\K,\delta\cotimes\id}
& 
& 
\\
{B\otimes\K,\epsilon\cotimes\id}
\ar[urr]^{\phi\otimes\id}
\ar[rrrr]^{\Xi}
&
&
&
&
{C\otimes\K,\eta\cotimes\id.} 
\ar[ull]_{\theta\otimes\id}
}
\]
Since $\epsilon$ and $\eta$ are assumed to be maximal, the outer two
vertical arrows are equivariant isomorphisms.  
Since $\epsilon$ and $\eta$ are
maximalizations of $\delta$, the upper two horizontal arrows are also
equivariant isomorphisms.  
Thus, an isomorphism $\Xi\colon B\otimes\K\to C\otimes\K$ can be
defined so that the outer rectangle commutes equivariantly. 
The two inner quadrilaterals commute equivariantly by straightforward
calculation, and therefore the lower triangle does as well. 

Now, \eqref{Phi-eq} (applied to $B$) says that $\Phi_B(p_B) = 1_B\otimes
P$, where $p_B = (k_{C(G)}^B\times k_G^B)(q)$,
and moreover,
\[
p_A 
= (k_{C(G)}^A\times k_G^A)(q)
= (\phi\times G\times G)\circ(k_{C(G)}^B\times k_G^B)(q)
= (\phi\times G\times G)(p_B).
\]
Combining this with similar calculations for $\Phi_C$ and $\theta$, 
we can conclude that $\Xi(1_B\otimes P) = 1_C\otimes P$.  
This, and the equivariant commutativity of the lower triangle, shows that
$\chi = \Xi|_{1\otimes P}$ is an $\epsilon - \eta$ equivariant isomorphism
of $B$ onto $C$ such that $\theta\circ\chi = \phi$.
\end{proof}

We remark that 
although any two maximalizations of $\delta$ are isomorphic, we don't
yet know how to define a \emph{canonical} maximalization. 
The construction of the maximalization $(A^m,G,\delta^m)$ in the proof
of \thmref{max thm} depended on the arbitrary choice of a rank-one
projection in $\K(L^2(G))$. 

Nonetheless, we feel that the basic approach of the proof will be
useful in greater generality; specifically, for Hopf $C^*$-algebras.
We have chosen not to treat Hopf algebras here because amenability
considerations are likely to muddy the waters significantly.

Also, we point out that 
while our results allow us to define an ``intermediate''
crossed product ``$\times_\mu$'' for \emph{dual actions} such that
$A\times_{\delta}G\times_{\what\delta,\mu}G$ is always naturally
isomorphic to $A\otimes\K$, 
this construction is not well-defined on isomorphism classes of dual
actions.
For instance, let $(B,G,\beta)$ be any action
for which the regular representation is not faithful, and let
$\hat\beta_r$ be the dual coaction on the reduced crossed product
$B\times_{\beta,r}G$.  Then 
\[
(B\times_\beta G\times_{\hat\beta}G,G,\hat{\hat\beta})
\cong
(B\times_{\beta,r}G\times_{\hat\beta_r}G,G,\hat{\hat\beta_r})
\]
because $\hat\beta_r$ is a normalization of $\hat\beta$ 
(\cite[Corollary~3.4 and Proposition~3.8]{QuiFR}). 
However, 
\[
(B\times_\beta G\times_{\hat\beta}G)\times_{\hat{\hat\beta},\mu}G
\cong (B\times_\beta G)\otimes\K,
\]
while
\[
(B\times_{\beta,r}G\times_{\hat\beta_r}G)\times_{\hat{\hat\beta_r},\mu}G
\cong (B\times_{\beta,r}G)\otimes\K.
\]
So $\hat{\hat\beta}$ and $\hat{\hat\beta_r}$ are isomorphic actions with
nonisomorphic intermediate crossed products.


\section{Discrete coactions}
\label{discrete}

For a coaction $\delta$ of a \emph{discrete} group $G$, 
we can use the results of \cite{EQ-FD} to show that
the dual coaction $\delta^m$ of $G$ on the maximal cross-sectional 
algebra
$C^*(\A)$ of the corresponding Fell bundle $\A$ is a
maximalization of $\delta$, and it is the {unique} maximalization 
with the same underlying Fell bundle. 

In preparation for this, we prove the following technical lemma, an  
easy modification of \cite[Lemma 2.5]{SieME},  
to streamline the task of verifying that a linear map is a  
right-Hilbert bimodule homomorphism.  The essential idea is that a  
linear map between Hilbert modules which preserves inner products is  
automatically a module homomorphism.  Let $C$ and $D$ be  
$C^*$-algebras, let $Z$ be a (right) Hilbert $D$-module, and suppose  
$C$ is represented by adjointable operators on $Z$.  If $Z$ is full 
as  
a Hilbert $D$-module and the action of $C$ on $Z$ is nondegenerate, 
we  
say $Z$ is a \emph{right-Hilbert $C - D$ bimodule}.  We use the  
notation ${}_CZ_D$ to indicate the coefficient algebras.
If $C_0\subseteq C$ and $D_0\subseteq D$ are dense  
$^*$-subalgebras and $Z_0$ is a dense linear subspace of $Z$ such that $C_0  
Z_0 \cup Z_0 D_0 \subseteq Z_0$ and $\< Z_0, Z_0 \>\rip D \subseteq  
D_0$, we say ${}_CZ_D$ is the \emph{completion} of the  
\emph{right-pre-Hilbert bimodule} ${}_{C_0}(Z_0)_{D_0}$. 
 
\begin{lem} 
\label{module hom} 
Suppose ${}_CZ_D$ and ${}_EW_F$ are right-Hilbert bimodules such 
that  
${}_CZ_D$ is the completion of a right-pre-Hilbert bimodule  
${}_{C_0}(Z_0)_{D_0}$, and suppose we are given homomorphisms $\phi 
\:  
C \to E$ and $\psi \: D \to F$ and a linear map $\theta \: Z_0\to W$  
with dense range such that for all $c \in C_0$ and $z,w \in Z_0$ we  
have 
\begin{enumerate} 
\item 
$\theta(cz) = \phi(c) \theta(z)$, and 
 
\item 
$\< \theta(z),\theta(w) \>\rip{F} = \psi\bigl(\< z,w \>\rip{D}\bigr)$. 
\end{enumerate} 
Then $\theta$ extends uniquely to a right-Hilbert bimodule 
homomorphism  
of ${}_CZ_D$ onto ${}_EW_F$.  Moreover, if $Z$ and $W$ are actually 
imprimitivity bimodules, this extension of $\theta$ is an 
imprimitivity  
bimodule homomorphism. 
\end{lem} 
 
\begin{proof} 
The argument of \cite[Lemma 2.5]{SieME}
shows that $\theta$ is bounded and preserves 
the right module actions, which implies that $\theta$ extends uniquely to 
a right-Hilbert bimodule homomorphism from ${}_CZ_D$ to 
${}_EW_F$. 
Since the range of $\theta$ is dense and Hilbert module homomorphisms 
have closed 
range, the extension of $\theta$ is onto, giving the first assertion. 
Moreover, when $Z$ and $W$ are actually imprimitivity bimodules, the 
density of $\theta(Z_0)$ together with the compatibility of the left 
and 
right inner products in both $Z$ and $W$ imply that 
$\theta$ preserves the left inner products. 
\end{proof} 
 
\begin{prop} 
\label{max char}
Let $\delta\:A\to A\otimes C^*(G)$ be a coaction of a
discrete group $G$ and let $\A$ be the corresponding Fell bundle
over $G$. Then $\delta$ is maximal if and only if $A=C^*(\A)$.
In general, the dual coaction $\delta^m\:C^*(\A)\to C^*(\A)\otimes C^*(G)$ is
a maximalization of $\delta$.
\end{prop}

\begin{proof}
Let $X$ be the 
$C^*(\c A)\times_{\delta^m}G\times_{\what{\delta^m}}G 
- C^*(\c A)$ imprimitivity bimodule of \cite[Theorem~3.1]{EQ-FD}
(for the case $H = G$), which is a completion of the sectional algebra
$\Gamma_c(\c A\times G)$ of the Fell bundle $\c A\times G$.
Using the isomorphism $C^*(\c
A)\times_{\delta^m}G \cong A\times_\delta G$ of
\cite[Lemma~2.1]{EQ-IC},  
$X$ becomes an $A\times_\delta G\times_{\deltahat}G - C^*(\c A)$ 
imprimitivity
bimodule, with left action and right inner-product given, for
$(a_r,s,t)\in \Gamma_c(\A\times G\times G)\subseteq
A\times_\delta G\times_{\deltahat}G$ and $(a_s,t),(a_u,v)\in
\Gamma_c(\A\times G)\subseteq X$, by 
\begin{gather*}
(a_r,s,t)\cdot(a_u,v) = (a_ra_u,vt^{-1})\nonzero{st=uv},\ \textup{and}\\
\<(a_s,t),(a_u,v)\>\rip{C^*(\A)} = a_s^*a_u \nonzero{st=uv}.
\end{gather*}

Let $\psi$ denote the unique homomorphism of $C^*(\A)$ onto $A$ which
extends the identity map on $\Gamma_c(\A)$.  We must show that the
canonical surjection $\Phi\colon A\times_\delta G\times_{\deltahat}G
\to A\otimes\K$ is an isomorphism if and only if $\psi$ is; for this,
it suffices to construct an imprimitivity bimodule homomorphism of $X$
onto the $A\otimes\K - A$ imprimitivity bimodule
$A\otimes\ell^2(G)$ with coefficient homomorphisms $\psi$ and $\Phi$.

For $a_s\in\A$ and $t\in G$, define 
\[ 
\Theta(a_s,t) = a_s \otimes \chi_{st}. 
\] 
Then $\Theta$ extends uniquely to a linear map from $\Gamma_c(\c A  
\times G)$ to $A \otimes \ell^2(G)$ with dense range.  
Since $\Phi(a_s, t, r) = a_s \otimes \lambda_s M_{\chi_t} \rho_r$,
a routine computation with generators shows that 
\[ 
\Theta(c \cdot z) = \Phi(c) \cdot \Theta(z) 
\midtext{and} 
\< \Theta(z),\Theta(w) \>\rip{A} 
= \psi\bigl( \< z,w \>\rip{C^*(\c A)} \bigr)
\] 
for $c \in \Gamma_c(\c A  \times G \times G)\subseteq A\times G\times
G$ and $z, w \in \Gamma_c(\c A\times G)\subseteq X$.
Thus by \lemref{module hom} the triple $(\Phi,\Theta, \psi)$
extends uniquely to an  
imprimitivity bimodule homomorphism, as desired.

The above arguments show that $\delta$ is maximal if and only if 
$A=C^*(\A)$. The remaining part follows from \cite[Lemma 
2.1]{EQ-IC}.
\end{proof} 



\providecommand{\bysame}{\leavevmode\hbox to3em{\hrulefill}\thinspace}

\end{document}